\documentclass{article}

\usepackage{arxiv}

\usepackage[utf8]{inputenc} 
\usepackage[T1]{fontenc}    
\usepackage{hyperref}       
\usepackage{url}            
\usepackage{booktabs}       
\usepackage{amsfonts}       
\usepackage{nicefrac}       
\usepackage{microtype}      
\usepackage{lipsum}
\usepackage{graphicx}
\graphicspath{ {./images/} }

\usepackage{amsmath,amssymb}   
\title{Epsilon-Optimal Policies for Average-Cost Separable MDPs with Perturbations}

\author{
 Dhairya Kantawala \\
  Indian Institute of Technology Bombay\\
  \texttt{dhairya@iitb.ac.in}
}

\begin{document}
\maketitle
\begin{abstract}
We study a class of infinite-horizon average-cost Markov Decision Processes (MDPs) whose reward and transition structures are nearly separable. For the totally separable baseline (i.e., with no perturbation), we derive an explicit stationary decision rule that is exactly average-optimal. We then show that under an $\varepsilon$-perturbation of the separable structure, this policy remains $\varepsilon$-optimal—that is, the loss in the average reward is of order $\mathcal{O}(\varepsilon)$.
\end{abstract}


\section{Introduction}

Markov Decision Processes (MDPs) provide a general framework for modeling sequential decision-making under uncertainty. An MDP consists of a finite set of states $\mathcal{S}$, a finite set of actions $\mathcal{A}$, a transition kernel $P(\cdot \mid s,a)$ specifying the probability of moving to each next state given the current state $s$ and action $a$, and a reward function $r(s,a)$ that determines the immediate reward obtained. The goal is to find a stationary policy $\pi : \mathcal{S} \to \mathcal{A}$ that optimizes a long-run performance criterion.

In this work, we focus on the \emph{average-cost} criterion, which measures the steady-state performance of a policy over an infinite horizon. Specifically, for a given stationary policy $\pi$, the average reward is defined as
\begin{equation}
    g_\pi = \lim_{T \to \infty} \frac{1}{T} \mathbb{E}_\pi \left[ \sum_{t=0}^{T-1} r(s_t, a_t) \right],
\end{equation}
where the expectation is taken with respect to the Markov chain induced by $\pi$ and $P$. Under mild ergodicity conditions, this limit exists and is independent of the initial state \cite{Puterman1994}. The standard objective is to find a stationary policy $\pi^*$ that maximizes $g_\pi$.

While the general theory of average-cost MDPs is well established \cite{Puterman1994,Bertsekas2012}, explicit closed-form expressions for the optimal policy or bias function are rare and usually require strong structural assumptions. This paper studies one such structured class, namely \emph{separable MDPs}, in which the reward and transition dynamics decompose additively across the state and action components. In the \emph{totally separable} case, the system takes the form
\begin{equation}
    r(s,a) = r_S(s) + r_A(a), \qquad 
    P(\cdot \mid s,a) = P_A(\cdot \mid a),
\end{equation}
which removes any direct state–action coupling from the dynamics. This decomposition leads to a tractable characterization of the average-optimal policy.

We first show that in the totally separable setting, the average-optimal stationary policy is explicit and constant across all states—the same action is optimal regardless of the starting state. We then extend the analysis to the \emph{nearly separable} or $\varepsilon$–separable case, where both the reward and transition kernels include small perturbations:
\begin{equation}
    r(s,a) = r_S(s) + r_A(a) + \varepsilon\,r_\varepsilon(s,a), \qquad 
    P(\cdot \mid s,a) = P_A(\cdot \mid a) + \varepsilon\,Q(\cdot \mid s,a).
\end{equation}
We prove that the explicit constant policy derived for the totally separable MDP remains $\varepsilon$–optimal under these perturbations, with the loss in average reward scaling as $\mathcal{O}(\varepsilon)$. This establishes a robustness guarantee for separable models and quantifies how weak state–action coupling affects long-run performance.

Our contributions are threefold:
\begin{enumerate}
    \item We derive an explicit closed-form expression for the average-optimal policy in totally separable MDPs.
    \item We establish a first-order perturbation result showing that this policy remains $\varepsilon$–optimal in the nearly separable regime.
\end{enumerate}

Throughout this work, we assume that under any stationary policy $\pi$, the induced Markov chain on $\mathcal{S}$ is \emph{unichain} and \emph{irreducible}. These standard assumptions ensure the existence and uniqueness (up to an additive constant) of the average-cost bias function and the steady-state distribution associated with each policy.

\section{Model Setup and Notation}

We consider a finite Markov Decision Process (MDP) defined by the tuple
\[
(\mathcal{S}, \mathcal{A}, P, r),
\]
where:
\begin{itemize}
    \item $\mathcal{S} = \{1, 2, \ldots, N\}$ is the finite set of states,
    \item $\mathcal{A} = \{1, 2, \ldots, M\}$ is the finite set of actions,
    \item $P(\cdot \mid s,a)$ is the transition kernel specifying the probability distribution over next states given $(s,a)$,
    \item $r(s,a)$ is the one-stage reward associated with taking action $a$ in state $s$.
\end{itemize}

A \emph{stationary policy} $\pi : \mathcal{S} \to \mathcal{A}$ assigns to each state a fixed action.  
Under a given stationary policy $\pi$, the controlled process $\{s_t\}_{t \ge 0}$ evolves as a time-homogeneous Markov chain with transition matrix
\[
P_\pi(s' \mid s) = P(s' \mid s, \pi(s)).
\]
The corresponding long-run average reward is defined as
\begin{equation}
    g_\pi = \lim_{T \to \infty} \frac{1}{T} \mathbb{E}_\pi \left[ \sum_{t=0}^{T-1} r(s_t, \pi(s_t)) \right],
\end{equation}
provided the limit exists. The objective is to find a stationary policy $\pi^*$ that maximizes $g_\pi$.

\subsection{Separable and Nearly Separable Structure}

We focus on a structured class of MDPs in which both the reward and transition dynamics are \emph{nearly separable} across the state and action components.  
Specifically, we consider models of the form
\begin{equation}
    r(s,a) = r_S(s) + r_A(a) + \varepsilon\,r_\varepsilon(s,a), 
    \qquad 
    P(\cdot \mid s,a) = P_A(\cdot \mid a) + \varepsilon\,Q(\cdot \mid s,a),
    \label{eq:separable}
\end{equation}
where $\varepsilon \ge 0$ is a small perturbation parameter.  
The case $\varepsilon = 0$ corresponds to a \emph{totally separable} MDP, in which the state and action influence the reward and transition independently.

The term $r_\varepsilon(s,a)$ represents a bounded perturbation to the additive reward structure, and $Q(\cdot \mid s,a)$ represents a bounded perturbation to the transition kernel.  
Throughout this work, we assume that for each $a \in \mathcal{A}$, $P_A(\cdot \mid a)$ defines a valid transition probability distribution, and that for all $(s,a)$, the perturbation kernel satisfies
\[
\sum_{s' \in \mathcal{S}} Q(s' \mid s,a) = 0.
\]
This ensures that $P(\cdot \mid s,a)$ remains a stochastic matrix for all $\varepsilon$ sufficiently small.

\subsection{Assumptions}

We impose the following standing assumptions, which are standard in the theory of average-cost MDPs \cite{Puterman1994,Bertsekas2012,HernandezLerma1999}:
\begin{enumerate}
    \item The state and action spaces $\mathcal{S}$ and $\mathcal{A}$ are finite.
    \item The reward function $r(s,a)$ is bounded for all $(s,a)$.
    \item Under every stationary policy $\pi$, the induced Markov chain with transition matrix $P_\pi$ is \emph{unichain} and \emph{irreducible}.
\end{enumerate}

The unichain and irreducibility assumptions guarantee that each stationary policy $\pi$ induces a unique invariant distribution $\pi_\pi$ and that the average reward $g_\pi$ is independent of the initial state.  
These conditions are sufficient to ensure the existence of a well-defined bias function and the validity of the average-cost optimality equations that will be used in subsequent sections.

\section{Baseline: Totally Separable MDP}

We now consider the baseline case $\varepsilon = 0$ in \eqref{eq:separable}, referred to as the \textit{totally separable MDP}. In this case, the reward and transition kernel decompose completely as
\[
r(s,a) = r_S(s) + r_A(a), \qquad P(\cdot \mid s,a) = P_A(\cdot \mid a),
\]
so that the action affects only the transition kernel, and the reward separates additively into a state component and an action component.

\textbf{Theorem 1.} (Totally Separable MDP)  
For the totally separable MDP defined above, suppose that for each action $a \in \mathcal{A}$, the Markov chain with transition matrix $P_A(\cdot \mid a)$ is irreducible with invariant distribution $\pi_a$. Then, the average-optimal stationary policy is constant across all states, i.e.
\[
\pi^*(s) = a^*, \qquad \forall s \in \mathcal{S},
\]
where
\[
a^* = \arg\max_{a \in \mathcal{A}} \left[\, \pi_a^\top r_S + r_A(a) \,\right],
\]
and the corresponding optimal average reward is
\[
g^* = \max_{a \in \mathcal{A}} \left[\, \pi_a^\top r_S + r_A(a) \,\right].
\]
Moreover, the pair $(g^*, h^*)$ satisfies the average-cost optimality equation exactly, i.e. $B(g^*, h^*) = 0$.

\textit{Proof.}  

1. {\bf ACOE in the totally separable case.}  
For $\varepsilon = 0$, the average-cost optimality equation (ACOE) for the gain $g$ and bias $h$ is
\[
g + h(s) = \max_{a \in \mathcal{A}} \Big[\, r_S(s) + r_A(a) + \sum_{s'} P_A(s' \mid a)\, h(s') \,\Big].
\]
(We use the standard ACOE form for average-cost MDPs; existence of $g,h$ under each stationary policy follows from the irreducibility/unichain assumption; see \cite{Puterman1994,Bertsekas2012}.)

2. {\bf Separation of state and action dependence.}  
Observe that the right-hand side decomposes as
\[
r_S(s) + \underbrace{\Big( r_A(a) + \sum_{s'} P_A(s' \mid a)\, h(s') \Big)}_{\text{term independent of } s}.
\]
Hence the dependence on $s$ appears only in the additive term $r_S(s)$ and does not affect the maximization over $a$.

3. {\bf The maximizing action is state-independent.}  
Since the inner maximization
\[
\max_{a \in \mathcal{A}} \Big\{ r_A(a) + \sum_{s'} P_A(s' \mid a)\, h(s') \Big\}
\]
does not depend on $s$, the same action attains the maximum for every state. Denote this maximizing action by $a^*$. Therefore the optimal stationary policy is constant:
\[
\pi^*(s) = a^*, \quad \forall s \in \mathcal{S}.
\]

4. {\bf Average reward under a constant policy.}  
Fix any constant policy $\pi(s)\equiv a$. By irreducibility of $P_A(\cdot\mid a)$ the Markov chain induced by $a$ has a unique invariant distribution $\pi_a$ (existence and uniqueness of invariant distribution; see \cite{Puterman1994}). The long-run average reward under this policy is the steady-state expectation of the one-stage reward:
\[
g(a) \;=\; \pi_a^\top \big( r_S + r_A(a)\,\mathbf{1} \big)
\;=\; \pi_a^\top r_S + r_A(a),
\]
where $\mathbf{1}$ denotes the all-ones vector and we used $\pi_a^\top \mathbf{1} = 1$. (This equality is the standard relation between average reward and invariant distribution; see \cite{Puterman1994}.)

5. {\bf Identification of $a^*$ and $g^*$.}  
Maximizing $g(a)$ over $a\in\mathcal{A}$ yields
\[
a^* = \arg\max_{a\in\mathcal{A}} \big[\, \pi_a^\top r_S + r_A(a)\,\big], \qquad
g^* = \max_{a\in\mathcal{A}} \big[\, \pi_a^\top r_S + r_A(a)\,\big].
\]
(Thus the constant action that maximizes the steady-state expected reward is optimal among constant policies; combined with step 3 this is the globally optimal stationary policy.)

6. {\bf Existence of the bias function for $a^*$.}  
We now show that there exists a bias function $h^*$ (unique up to an additive constant) satisfying the Poisson/bias equation for $a^*$:
\[
g^* + h^*(s) = r_S(s) + r_A(a^*) + \sum_{s'} P_A(s' \mid a^*)\, h^*(s').
\]
Rewriting in vector form gives
\[
(I - P_A(a^*))\, h^* \;=\; r_S + r_A(a^*)\,\mathbf{1} - g^* \mathbf{1}.
\]
Solvability of this linear system requires the right-hand side to be orthogonal to the left nullspace of $(I-P_A(a^*))$, i.e. to $\pi_{a^*}^\top$. Indeed,
\[
\pi_{a^*}^\top \big( r_S + r_A(a^*)\,\mathbf{1} - g^* \mathbf{1} \big)
= \pi_{a^*}^\top r_S + r_A(a^*) - g^* = g(a^*) - g^* = 0,
\]
since $g^* = g(a^*)$ by definition. (Hence the compatibility condition holds.) By standard Poisson equation theory for irreducible chains there exists a solution $h^*$, unique up to adding a constant; see \cite{HernandezLerma1999,Puterman1994}.

7. {\bf Verification of the ACOE and conclusion.}  
With $a^*$ and $h^*$ as above we have, for every state $s$,
\[
g^* + h^*(s) = r_S(s) + r_A(a^*) + \sum_{s'} P_A(s' \mid a^*)\, h^*(s'),
\]
which matches the ACOE since the right-hand side is the maximized value over actions (the maximizer being $a^*$ independent of $s$). Therefore $(g^*,h^*)$ satisfies the ACOE exactly, and by definition $B(g^*,h^*)=0$.

\section{Perturbed Case: Nearly Separable MDP}

We now consider the nearly separable case with $\varepsilon > 0$ in \eqref{eq:separable}. Our goal is to show that the constant policy derived for the totally separable MDP remains nearly optimal under small perturbations.

\textbf{Theorem 2.} (Nearly Separable MDP)  
For the nearly separable MDP defined in \eqref{eq:separable}, let $\pi^*(s) = a^*$ be the constant optimal policy from Theorem 1, where $a^*$ maximizes $\pi_a^\top r_S + r_A(a)$. Assume the perturbation terms $r_\varepsilon(s,a)$ and $Q(\cdot \mid s,a)$ are bounded uniformly over all $(s,a)$. Then, there exist constants $C > 0$ and $\varepsilon_0 > 0$ such that for all $0 \le \varepsilon < \varepsilon_0$, the policy $\pi^*$ is $\mathcal{O}(\varepsilon)$-optimal in the perturbed MDP. Specifically,
\[
g_\varepsilon^* - g_{\pi^*,\varepsilon} \le C \varepsilon,
\]
where $g_\varepsilon^*$ is the optimal average reward in the $\varepsilon$-perturbed MDP, and $g_{\pi^*,\varepsilon}$ is the average reward under $\pi^*$ in the perturbed MDP.

\textit{Proof.}  

1. {\bf Perturbation bound for fixed policies.}  
For any fixed stationary policy $\pi : \mathcal{S} \to \mathcal{A}$, denote by $g_\pi^0$ and $g_\pi^\varepsilon$ the average rewards under $\pi$ in the unperturbed ($\varepsilon=0$) and perturbed models, respectively. Similarly, let $P_\pi^0$ and $P_\pi^\varepsilon = P_\pi^0 + \varepsilon Q_\pi$ be the corresponding transition matrices, where $Q_\pi(s' \mid s) = Q(s' \mid s, \pi(s))$, and let $r_\pi^0(s) = r_S(s) + r_A(\pi(s))$ and $r_\pi^\varepsilon(s) = r_\pi^0(s) + \varepsilon r_\varepsilon(s, \pi(s))$. Let $\pi^0$ and $\pi^\varepsilon$ be the unique invariant distributions of $P_\pi^0$ and $P_\pi^\varepsilon$.

Under the irreducibility assumption, for sufficiently small $\varepsilon$, the perturbed chain remains irreducible and unichain. By standard perturbation theory for Markov chains \cite{Meyer1980,Seneta1981}, the invariant distribution admits a regular expansion:
\[
\pi^\varepsilon = \pi^0 + \varepsilon \pi^1 + O(\varepsilon^2),
\]
where the first-order correction is given by
\[
\pi^1 = -\pi^0 Q_\pi (I - P_\pi^0)^\#,
\]
and $(I - P_\pi^0)^\#$ is the group inverse of $(I - P_\pi^0)$, which exists and is bounded under the unichain assumption.

2. {\bf Expansion of the average reward.}  
The perturbed average reward is
\[
g_\pi^\varepsilon = \pi^\varepsilon r_\pi^\varepsilon = (\pi^0 + \varepsilon \pi^1 + O(\varepsilon^2)) (r_\pi^0 + \varepsilon r_{\varepsilon,\pi}),
\]
where $r_{\varepsilon,\pi}(s) = r_\varepsilon(s, \pi(s))$. Expanding gives
\[
g_\pi^\varepsilon = \pi^0 r_\pi^0 + \varepsilon (\pi^0 r_{\varepsilon,\pi} + \pi^1 r_\pi^0) + O(\varepsilon^2).
\]
Since $g_\pi^0 = \pi^0 r_\pi^0$ and all terms (rewards, transition perturbations, and the group inverse) are bounded, it follows that
\[
|g_\pi^\varepsilon - g_\pi^0| \le C_\pi \varepsilon
\]
for some constant $C_\pi > 0$ depending on $\pi$, valid for $\varepsilon$ sufficiently small.

3. {\bf Uniform bound over policies.}  
Since the state and action spaces are finite, there are only finitely many stationary policies (at most $M^N$). Therefore, we can take
\[
C = \max_\pi C_\pi < \infty,
\]
yielding a uniform bound $|g_\pi^\varepsilon - g_\pi^0| \le C \varepsilon$ over all stationary policies $\pi$, for $0 \le \varepsilon < \varepsilon_0$ with $\varepsilon_0$ small enough to ensure all perturbed chains remain unichain.

4. {\bf Bound on the optimal average reward.}  
The optimal average reward in the perturbed MDP is
\[
g_\varepsilon^* = \max_\pi g_\pi^\varepsilon \le \max_\pi (g_\pi^0 + C \varepsilon) = \max_\pi g_\pi^0 + C \varepsilon = g^* + C \varepsilon,
\]
where $g^* = \max_\pi g_\pi^0$ is the optimal unperturbed average reward.

5. {\bf Bound on the performance of $\pi^*$.}  
For the specific policy $\pi^*(s) = a^*$, we have $g_{\pi^*}^0 = g^*$ by Theorem 1, and thus
\[
g_{\pi^*}^\varepsilon \ge g_{\pi^*}^0 - C \varepsilon = g^* - C \varepsilon.
\]

6. {\bf Optimality gap.}  
Combining the above inequalities yields
\[
g_\varepsilon^* - g_{\pi^*}^\varepsilon \le (g^* + C \varepsilon) - (g^* - C \varepsilon) = 2C \varepsilon.
\]
Therefore, $\pi^*$ is $\mathcal{O}(\varepsilon)$-optimal with constant $C' = 2C$. 

\section{Conclusion}

We derived explicit average-optimal stationary policies for totally separable MDPs and showed that they remain $\mathcal{O}(\varepsilon)$–optimal under small perturbations of the separable structure. 
This result formalizes the intuition that separability leads to robustness in long-run decision-making. 

Future work could extend these results to:
(i) discounted-cost and constrained MDPs, 
(ii) continuous-state separable systems, and 
(iii) learning-based or empirical estimation settings where the perturbation parameter $\varepsilon$ reflects statistical uncertainty rather than model structure.

\bibliographystyle{unsrt}  
\bibliography{references}  

\end{document}